\documentclass[a4paper, 12pt]{amsart}

\usepackage{amssymb}
\usepackage{amsthm}
\usepackage{amsmath}
\usepackage{mathrsfs}
\usepackage{comment}

\AtBeginDocument{%
   \def\MR#1{}
}

\theoremstyle{plain}
\newtheorem{thm}{Theorem}[section]		
\newtheorem{prop}[thm]{Proposition}
\newtheorem{cor}[thm]{Corollary}
\newtheorem{lem}[thm]{Lemma}

\theoremstyle{definition}
\newtheorem{df}{Definition}[section]

\theoremstyle{remark}
\newtheorem{rmk}{Remark}[section]
\newtheorem*{ac}{Acknowledgements}

\newcommand{\zz}{\mathbb{Z}}

\newcommand{\rr}{\mathbb{R}}

\DeclareMathOperator{\met}{Met}

\DeclareMathOperator{\ult}{UMet}

\newcommand{\opcpt}[1]{\widetilde{\alpha}#1}

\newcommand{\orddis}{M}

\newcommand{\eucdis}{E}

\DeclareMathOperator{\conv}{\mathcal{CC}}

\DeclareMathOperator{\clo}{\mathcal{C}}

\newcommand{\yosub}{\subset}

\newcommand{\yopow}[1]{P(#1)}

\makeatletter
\@addtoreset{equation}{section}

\makeatother

\begin{document}

\title[Extending proper metrics]
{
Extending proper metrics
}
\author[Yoshito Ishiki]
{Yoshito Ishiki}
\address[Yoshito Ishiki]
{\endgraf
Photonics Control Technology Team
\endgraf
RIKEN Center for Advanced Photonics
\endgraf
2-1 Hirasawa, Wako, Saitama 351-0198, Japan}
\email{yoshito.ishiki@riken.jp}

\date{\today}

\maketitle

\begin{abstract}
We first prove 
a version of 
Tietze-Urysohn's theorem for 
proper functions taking values in non-negative 
real 
numbers defined on  $\sigma$-compact locally compact Hausdorff  spaces. 
As an application, 
we prove an extension theorem of 
proper metrics, 
which states that 
if  $X$ is  a 
 $\sigma$-compact locally compact Hausdorff space, $A$ is a  closed subset of $X$, 
 and 
$d$ is  a 
proper metric on $A$ that generates the same topology of $A$, 
then 
 there exists a proper metric $D$ on $X$ such that
$D$ generates the same topology of $X$ and 
 $D|_{A^{2}}=d$. 
Moreover, if $A$ is a proper retract, 
we can choose $D$ so that 
$(A, d)$ is quasi-isometric  to  $(X, D)$. 
We also show analogues of the 
theorems explained above for ultrametrizable  spaces. 
\end{abstract}

\section{Introduction}

Tietze--Urysohn's theorem  states that 
every  continuous function on a closed subset of a normal space can be  extended to the whole space
as a continuous function. 
This theorem
has played an important role in topology and analysis 
(for example,  the existence of a partition of unity). 
There are many 
generalizations of 
Tietze--Urysohn's theorem 
 (see for instance  \cite{Du1951}, 
\cite{MR1967003}, \cite{MR77107} and \cite{MR1346246}). 

For a metrizable space $X$, 
we denote by $\met(X)$ the 
set of all metrics on $X$ generating the 
same topology of $X$. 
 Hausdorff's extension theorem 
 states that for every metrizable space $X$, 
 for every  closed subset $A$ of $X$, 
 and for every $d\in \met(A)$, 
 there exists a metric $D\in \met(X)$ such that 
 $D|_{A^{2}}=d$. 
 This theorem can be considered as 
 an analogue of 
 Tietze--Urysohn's theorem for metric spaces, 
 and some variants have been 
  investigated by some authors
 (see \cite{Bing1947}, \cite{MR3135687}, \cite{MR3090172}, and \cite{Ishiki2021ultra}). 
 
 Tietze--Urysohn's theorem  and 
 Hausdorff's extension theorem 
 are not only analogous, 
 but also logically connected with each other. 
 In fact, 
 according to \cite{Arens1952}, 
 Hausdorff's extension theorem can be proven using Dugundji's theorem (see \cite{Du1951}), which is 
an improvement of Tietze--Urysohn's theorem. 
For more discussion on connections between 
extensions of maps and metrics, 
we refer the readers to \cite{Hu2010}. 

In  \cite{Ishiki2021ultra}, 
the author proved an ultrametric analogue of 
 Hausdorff's  extension theorem using 
the method described above, namely, using 
the property that 
every continuous function on a closed subsets of  an 
ultrametrizable space can be extended to the whole space. 
Remark that 
since every non-empty closed set in an ultrametric  space is a retract of the whole space (see \cite[Theorem 2.9]{brodskiy2007dimension}), 
all extension  problems of continuous maps
defined on a closed subsets  of an ultrametric  space are  solved affirmatively.

Let $X$ and $Y$ be topological spaces. 
A  map
$f\colon X\to Y$
 is said to be 
\emph{proper} if 
for every compact  subset $K$ of $Y$, 
the inverse image  $f^{-1}(K)$ is compact. 
A metric $d$ on $X$ is said to be \emph{proper} if 
all bounded closed subsets of $(X, d)$ are compact. 
In this case, for a fixed point $p\in X$, the function defined by 
$x\mapsto d(p, x)$ is a proper map. 
These two concepts are the 
main subjects of this paper.

In the present  paper, we prove a new variant of
 Hausdorff's extension theorem
for proper metrics. 
The key idea is, 
as mentioned above, 
that 
an  extension theorem of continuous maps 
implies  an 
 extension theorem of metrics. 
We first show   Tietze--Urysohn's theorem
for proper functions (Theorem \ref{thm:TUproper}), 
and then, as an application,  we prove 
Hausdorff's extension theorem for proper metrics
(Theorem \ref{thm:extproper}). 

To  show  Tietze--Urysohn's theorem for proper functions,  we use  the 
so-called controlling 
Tietze--Urysohn's theorem 
(see  \cite{MR1346246} and
\cite{MR1967003}), 
which includes
not only an extension of a given function, but also  an extension of the zero set of the function.

Similarly, 
using the fact that every non-empty closed subset of 
an ultrametric space is a retract of the whole space, 
we prove an extension theorem of proper ultrametrics
(Theorem \ref{thm:extproperult}).

In this paper, 
we also prove an extension theorem of 
proper metrics
focusing on large scale structures of  metric spaces
(Theorem \ref{thm:extdenseproper})
using Michael's continuous selection theorems, 
which are also generalizations of Tietze--Urysohn's theorem. 
More precisely, we prove that 
for every $\sigma$-compact locally compact 
 space $X$, and for 
every closed subset $A$ of  $X$, 
if $A$ is a proper retract of $X$, 
then for every proper metric $d\in \met(X)$, 
there exists a proper metric $D\in \met(X)$
such that $(A, D)$ is quasi-isometric to $(X, D)$. 
We also prove an
ultrametric  version of this extension theorem
(Theorem \ref{thm:extultproper}).

\section{Extension of proper functions}
A main purpose of this section is to  prove 
Tietze--Urysohn's 
theorem for proper functions 
(Theorem \ref{thm:TUproper}). 

A topological space is said to be 
\emph{$\sigma$-compact} if 
it is a countable  union of 
compact subspaces.
A topological space is said to be \emph{locally compact} if 
every point in the space has a compact neighborhood. 

\begin{rmk}
As a consequence of Urysohn's metrization theorem 
(see \cite[Theorem 34.1]{MR3728284}), 
all  $\sigma$-compact locally compact Hausdorff 
spaces are metrizable. 
Indeed, they are second countable and regular. 
\end{rmk}

The following
theorem  is 
deduced from Yamazaki's theorem 
\cite[Corollary 2.1]{MR1967003} or 
Frantz's theorem \cite[Theorem 1]{MR1346246}.  
\begin{thm}\label{thm:zeroTU}
Let $X$ be a normal space, 
$A$  a closed subset of $X$, 
and 
 $Z$  a closed $G_{\delta}$ subset of $X$. 
Assume that  $f\colon A\to [0, 1]$ is  a continuous function
such that $Z\cap A=f^{-1}(0)$. 
Then there exists a continuous function 
$F\colon X\to [0, 1]$ satisfying that 
$F|_{A}=f$ and $F^{-1}(0)=Z$. 
\end{thm}

For a $\sigma$-compact locally compact Hausdorff space $X$, 
we put 
$\opcpt{X}=X\sqcup \{\infty\}$. 
We define  a topology on $\opcpt{X}$
by declaring  neighborhood systems 
of $\opcpt{X}$ 
as follows: 
If $p\in X$, 
then 
the neighborhood system of $p$ in  $\opcpt{X}$ is 
the family  of all subsets $V$ of $X$ such that 
$V\setminus \{\infty\}$ is a neighborhood of $p$ in $X$, 
and  
the neighborhood system of $\infty$ in $\opcpt{X}$
is 
the set of all subsets  $V$ of $\opcpt{X}$ 
satisfying  that
$\opcpt{X}\setminus V$ is a 
relatively 
compact subset of $X$. 
In what follows, we always consider that 
$\opcpt{X}$ is equipped with 
this topology. 
If $X$ is non-compact, 
then $\opcpt{X}$ coincides  with the 
one-point compactification of $X$. 
If $X$ is compact, 
then 
$\opcpt{X}$ is nothing but the topological
direct  sum of $X$ and 
the point $\infty$. 
Remark that if $X=\emptyset$, 
then $\opcpt{X}=\{\infty\}$. 

Let $X$ and 
$Y$ be  $\sigma$-compact locally compact Hausdorff spaces. 
For a map $f\colon X\to Y$, 
we define an induced  map  $\opcpt{f}\colon \opcpt{X}\to \opcpt{Y}$ 
by 
$\opcpt{f}|_{X}=f$ and 
$\opcpt{f}(\infty)=\infty$. 
\begin{prop}\label{prop:alphaext}
Let $X$ and $Y$ be  $\sigma$-compact locally compact Hausdorff spaces. 
Then the following statements hold: 
\begin{enumerate}
\item\label{item:alpha1}
 For every proper map $f\colon X\to Y$, 
the map  $\opcpt{f}\colon \opcpt{X}\to \opcpt{Y}$ is continuous. 
\item\label{item:alpha2} 
If a continuous map $F\colon \opcpt{X}\to \opcpt{Y}$
satisfies $F^{-1}(\infty)=\{\infty\}$, 
then the restriction $F|_{X}\colon X\to Y$ is 
proper. 
\end{enumerate}
\end{prop}
\begin{proof}
We first prove \eqref{item:alpha1}. 
Let $A$ be a closed subset of $\opcpt{Y}$. 
Then $A$ is compact and it is  contained in $Y$, 
 or $A=B\cup\{\infty\}$ for some closed subset 
 $B$  of $Y$.
 In any case, the inverse image 
 $(\opcpt{f})^{-1}(A)$ is closed. 
Thus  $\opcpt{f}$ is continuous. 
  
 To prove \eqref{item:alpha2}, 
 we
 take an arbitrary compact subset $K$ of $Y$. 
 Since $\infty\not\in K$, 
 we have $F^{-1}(\infty)\cap F^{-1}(K)=\emptyset$. 
By  $F^{-1}(\infty)=\{\infty\}$, 
we obtain $\infty\not\in F^{-1}(K)$. 
This means that $F^{-1}(K)$ is compact in $X$. 
Thus $F|_{X}$ is proper.  
\end{proof}

\begin{thm}\label{thm:TUproper}
Let $X$ be a 
$\sigma$-compact locally compact 
Hausdorff space, 
and  $A$  a  
closed subset of $X$. 
If 
$f\colon A\to  [0, \infty)$ is  a continuous proper function, 
then there exists 
a continuous proper function $F\colon X\to [0, \infty)$
such that $F|_{A}=f$. 
\end{thm}
\begin{proof}
Notice that  $\opcpt{[0, \infty)}$ is homeomorphic to $[0, \infty]$. 
According to  \eqref{item:alpha1} in Proposition \ref{prop:alphaext}, 
the map $\opcpt{f}\colon \opcpt{A}\to [0, \infty]$ is 
continous. 
The space $\opcpt{A}$ can be considered as a
closed subset of $\opcpt{X}$. 
Since $X$ is $\sigma$-compact, 
the singleton $\{\infty\}$ is 
a
closed $G_{\delta}$ set in $\opcpt{X}$. 
The space   $[0, \infty]$
 is homeomorphic to $[0, 1]$. 
 Since $\opcpt{X}$ is compact  and 
  Hausdorff, 
 it is normal. 
Thus, due to  Theorem \ref{thm:zeroTU}, 
there exists a continuous 
map $h\colon \opcpt{X}\to [0, \infty]$ such that 
$h|_{\opcpt{A}}=\opcpt{f}$ and $h^{-1}(\infty)=\{\infty\}$. 
By \eqref{item:alpha2} in Proposition \ref{prop:alphaext}, 
the function $F=h|_{X}\colon X\to [0, \infty)$
 is proper and satisfies $F|_{A}=f$. 
 This finishes the proof of 
 Theorem \ref{thm:TUproper}. 
\end{proof}

\begin{rmk}
In Theorem \ref{thm:TUproper}, 
it is important  that the target space is 
$[0, \infty)$. 
In general, a  proper function $f\colon A\to \rr$
can not be extended to the ambient space as a proper function. 
For example, 
if 
we define a map $f\colon \zz\to \rr$
 by $f(n)=(-1)^{n}\cdot n$, then $f$ is proper. 
However, for any continuous extension $F\colon \rr\to \rr$ of $f$, 
the set $F^{-1}(0)$ is non-compact  by 
the intermediate value theorem. 
\end{rmk}

The following proposition is well-known. 
However, 
for the sake of self-containedness, 
we provide a proof. 
\begin{prop}\label{prop:sigmasigma}
A Hausdorff space is 
$\sigma$-compact and locally compact if and only if 
 there exists a continuous proper function 
$f\colon X\to [0, \infty)$. 
\end{prop}
\begin{proof}
We first assume that $X$ is $\sigma$-compact and 
locally compact. 
Applying Theorem \ref{thm:TUproper} to 
$A=\emptyset$ and the empty map from $\emptyset$ into 
$[0, \infty)$, 
we obtain a proper function from $X$ into $[0, \infty)$. 

Next assume that there exists a continuous proper 
function $f\colon X\to [0, \infty)$. 
By  $X=\bigcup_{i=0}^{\infty}f^{-1}([0, i])$, 
the space $X$ is $\sigma$-compact. 
Since 
$X=\bigcup_{i=0}^{\infty}f^{-1}([0, i))$ and 
each $f^{-1}([0, i))$ is open and relatively compact, 
the space $X$ is locally compact. 
\end{proof}

\section{Extension of proper metrics}
In this section, 
we shall prove two extension theorems of proper metrics and 
ultrametrics 
(Theorems \ref{thm:extproper}
and \ref{thm:extproperult}).

A metric $d$ on $X$ is said to be 
\emph{ultrametric} if it satisfies 
$d(x, y)\le d(x, z)\lor d(z, y)$ for all $x, y, z\in X$, 
where $\lor$ is the maximum operator on $\rr$. 
A topological space is said to be 
\emph{metrizable} (resp.~\emph{ultrametrizable})
 if there exists a metric (resp.~ultrametric)
that generates the same topology of the space. 
Let $X$ be a metrizable space, 
and $S$  a subset of $[0, \infty)$ with 
$0\in S$. 
We denote by $\met(X; S)$
(resp.~$\ult(X; S)$)
the set of all metrics (resp.~ultrametrics) that 
generate the same topology of $X$ taking 
values in $S$. 
We often write 
$\met(X)=\met(X; [0, \infty))$.

A topological space 
$X$ is said to be \emph{ultranormal}
if for every pair $A$ and $B$ of disjoint closed subsets of $X$, 
there exists a clopen set $V$ such that 
$A\yosub V$
and $V\cap B=\emptyset$. 
Note that 
a topological space $X$ is 
ultrametrizable if and only if 
it is metrizable and ultranormal
(see \cite[Theorem II]{MR80905}). 
For a topological space $X$, 
a pair of subsets  $A$ and $B$ of $X$ is  
said to be \emph{completely separated} if 
there exists a continuous function 
$f\colon X\to [0, 1]$ such that $f^{-1}(0)=A$ and 
$f^{-1}(1)=B$. 
A topological space 
$X$ 
is \emph{strongly $0$-dimensional} if 
$X$ is completely regular and any two completely separated subsets of $X$ are separated by a clopen subset of 
$X$. 
Remark that the class of ultranormal spaces coincides with  the class of strongly $0$-dimensional normal spaces.
In particular, 
a metrizable space  is ultranormal if and only if it is strongly $0$-dimensional.

The next is  Hausdorff's extension theorem \cite{Ha1930}
(see also \cite{Bing1947} and  \cite{MR321026}). 
\begin{thm}\label{thm:Hausdorff}
Let $X$ be  a metrizable space, 
and  $A$ 
a closed subset $A$ of $X$. 
If $d\in \met(A)$, 
then there exists $D\in \met(X)$
such that $D|_{A^{2}}=d$. 
\end{thm}

A subset $S$ of $[0, \infty)$ is said to be 
\emph{characteristic}
 if $0\in S$ and 
if  for all $r\in (0, \infty)$, 
there exists $s\in S\setminus \{0\}$ with 
$s\le r$. 

We next  explain the author's extension theorem  of ultrametrics \cite[Theorem 1.2]{Ishiki2021ultra}, 
which 
is an analogue of  Hausdorff's  extension theorem for ultrametrics: 
\begin{thm}\label{thm:Ishiki}
Let $S$ be a characteristic subset of $[0, \infty)$. 
Let $X$ be  an ultrametrizable space, 
and 
$A$ 
 a closed subset of $X$. 
If 
$d\in \ult(A; S)$, 
then 
there exists 
$D\in \ult(X; S)$
such that $D|_{A^{2}}=d$. 
\end{thm}

The following  proposition can be 
considered as  a
$0$-dimensional analogue of 
Proposition \ref{prop:sigmasigma}. 
\begin{prop}\label{prop:intoS}
Let $S$ be an unbounded  subset of $[0, \infty)$, 
and 
$X$  an ultranormal 
$\sigma$-compact  locally compact Hausdorff space. 
Then 
 there exists a continuous proper function 
$f\colon X\to S$. 
\end{prop}
\begin{proof}
Let $\{U_{i}\}_{i\in I}$ be an 
open covering of $X$ consisting of 
relatively compact subsets. 
Since $X$ is paracompact and 
ultranormal, 
using 
\cite[Corollary 1.4]{ellis1970extending}, 
we obtain an open covering $\{V_{j}\}_{j\in J}$
of $X$
refining $\{U_{i}\}_{i\in I}$
such that 
$V_{j}\cap V_{j^{\prime}}= \emptyset$ if $j\neq j^{\prime}$. 
In this case, each $V_{j}$ is clopen and compact. 
Since $X$ is $\sigma$-compact, 
the set $J$ is at most  countable. 
We may assume that 
$J\yosub \zz_{\ge 0}$. 
Take a strictly increasing sequence 
$\{a_{j}\}_{j\in \zz_{\ge 0}}$
taking values 
in $S$ such that $\lim_{j\to \infty}a_{j}=\infty$. 
We define a map $f\colon X\to S$ by 
$f(x)=a_{j}$ if $x\in V_{j}$. 
From the fact that $\{V_{j}\}_{j\in J}$ is a mutually disjoint clopen covering of $X$, 
it 
follows that 
the map $f$ is continuous. 
Since each $V_{j}$ is compact, 
we conclude  that $f$ is proper. 
\end{proof}

Recall that  
the symbol $\lor$ stands for the 
maximum operator on $\rr$. 
Namely, $x\lor y=\max\{x, y\}$. 
\begin{df}\label{df:kappasp}
Let $S$ be a subset of $[0, \infty)$ with $0\in S$. 
We define an ultrametric 
$\orddis_{S}$ by 
\[
\orddis_{S}(x, y)
=
\begin{cases}
0 & \text{if $x=y$;}\\
x\lor y & \text{if $x\neq y$. }
\end{cases}
\]
\end{df}

\begin{rmk}
The construction of the metric  $\orddis_{S}$ was given by 
 Delhomm\'{e}--Laflamme--Pouzet--Sauer
 \cite[Proposition 2]{MR2435142}, 
which also can be found in \cite{Ishiki2021ultra} and
\cite{MR2854677}. 
\end{rmk}

Let $(X, d)$ be  a metric space, 
 $x\in X$, 
and 
$\epsilon\in (0, \infty)$. 
We denote by 
$U(x, \epsilon; d)$ (resp.~$B(x, \epsilon; d)$)
the open (resp.~closed) ball centered at 
$x$ with radius $\epsilon$. 

A subset $S$ of $[0, \infty)$ is said to be 
\emph{sporadic} 
if 
there exists a sequence $\{s_{n}\}_{n\in \zz}$ such that 
$S=\{0\}\cup \{\,s_{n}\mid n\in \zz \, \}$, 
$\lim_{n\to -\infty} s_{n}=0$, $\lim_{n\to \infty}s_{n}=\infty$
and $s_{i}<s_{i+1}$  for all $i\in \zz$.  
Note that every sporadic subset of $[0, \infty)$ is unbounded
 and characteristic 
in $[0, \infty)$, 
 and every unbounded characteristic 
subset of $[0, \infty)$ contains a 
sporadic subset. 

\begin{lem}\label{lem:sametop}
Let $S$ be a sporadic subset of $[0, \infty)$. 
Then the Euclidean topology on $S$ 
coincides  with 
that induced from $M_{S}$. 
\end{lem}
\begin{proof}
For all $x\in S\setminus \{0\}$, 
we have $U(x, x; \orddis_{S})=\{x\}$ and 
$U(0, x; \orddis_{S})=S\cap [0, x)$. 
This proves the lemma. 
\end{proof}

\begin{df}\label{df:mfmetric}
Let $X$ be a topological 
space, and $f\colon X\to \rr$  a continuous
 map. 
 We define a pseudo-metric 
 $\eucdis[f]$ on $X$
  by $\eucdis[f](x, y)=|f(x)-f(y)|$. 
Let $S$ be a subset of $[0, \infty)$. 
Let $f\colon X\to S$ be a continuous map. 
We also  define 
a pseudo-metric 
$\orddis_{S}[f]$ on $X$ by 
$\orddis_{S}[f](x, y)=\orddis_{S}(f(x), f(y))$. 
\end{df}
\begin{df}
Let $X$ be a set and 
 $d, e\colon X^{2}\to \rr$ 
 be arbitrary  maps. 
We define $d\lor e\colon X^{2}\to \rr$
 by $(d\lor e)(x, y)=d(x, y)\lor e(x, y)$. 
 Notice that if $d$ is a metric on $X$ and $e$ is a pseudo-metric on 
 $X$, then $d\lor e$ is a metric on $X$. 
\end{df}

Note that a metric $d$ on $X$ is proper if and only if 
all closed balls of $(X, d)$ are compact. 
\begin{lem}\label{lem:dormf}
Let $X$ be a metrizable space. 
Let $f\colon X\to [0, \infty)$ be a 
continuous
proper  function, 
and $d\in \met(X)$. 
Then the map  
$d\lor \eucdis[f]$ is  
a proper  metric in $\met(X)$. 
\end{lem}
\begin{proof}
Since $f$  is continuous, 
the map $\eucdis[f]\colon X^{2}\to [0, \infty)$
is also continuous. 
Then the assumption 
$d\in \met(X)$ yields 
$d\lor \eucdis[f]\in \met(X)$. 
For all $r\in (0, \infty)$ and $p\in X$,  we notice that 
$B(p, r; d\lor \eucdis[f])\yosub f^{-1}([f(p)-r, f(p)+r])$. 
Since $f$ is proper, the set $B(p, r; d\lor \eucdis[f])$ is compact. 
Thus, we conclude that $d\lor \eucdis[f]$ is a proper metric. 
\end{proof}

\begin{lem}\label{lem:dormfult}
Let $S$ be an unbounded characteristic subset  of $[0, \infty)$, 
and 
$T$  a sporadic subset of $[0, \infty)$ with 
$T\yosub S$. 
Let $X$ be an ultrametrizable space. 
Assume that $f\colon X\to T$ is  a 
continuous
proper  function, 
and $d\in \ult(X; S)$. 
Then the map  
$d\lor M_{T}[f]$ is 
a proper  metric in $\ult(X; S)$. 
\end{lem}
\begin{proof}
Lemma \ref{lem:sametop}
implies that $\orddis_{T}[f]\colon X^{2}\to T$ is 
continuous. Thus, by $d\in \ult(X; S)$, and by $T\yosub S$, 
we have $d\lor M_{T}[f]\in \ult(X; S)$. 
For all $r\in (0, \infty)$ and $p\in X$,  we obtain  
$B(p, r; d\lor \orddis_{T} [f])\yosub f^{-1}([0, r] \cup \{f(p)\})$. 
Since $f$ is proper, the set $B(p, r; d\lor \orddis_{T} [f])$ is compact. 
Therefore $d\lor M_{T}[f]$ is a proper metric. 
This completes the proof. 
\end{proof}

Lemma \ref{lem:dormf} gives  
a new proof of 
the following well-known corollary: 
\begin{cor}\label{cor:extpropermet}
Let $X$ be a $\sigma$-compact locally compact
Hausdorff  space. 
Then 
there exists a proper metric in 
$\met(X)$. 
In particular, 
the space  $X$ is completely metrizable. 
\end{cor}
\begin{proof}
Take $d\in \met(X)$ and take a proper continuous 
function $f\colon X\to [0, \infty)$ 
(see Proposition \ref{prop:sigmasigma}). 
By Lemma \ref{lem:dormf}, 
the map  $d\lor \eucdis[f]$ is a proper metric in $\met(X)$. 
This proves the first part of the corollary. 
The latter part follows from the fact that 
every proper metric is complete. 
\end{proof}

\begin{cor}\label{cor:extpropermetult}
Let $S$ be an unbounded characteristic subset of 
$[0, \infty)$, 
and 
$X$  an ultranormal 
$\sigma$-compact locally compact
Hausdorff space. 
Then there exists a proper metric  in 
 $\ult(X; S)$.
 In particular, 
 the space $X$ is completely ultrametrizable. 
\end{cor}
\begin{proof}
Since $S$ is unbounded and characteristic, 
there exists a sporadic set of $T$ such that $T\yosub S$. 
According to Proposition \ref{prop:intoS}, 
there exists a continuous proper function $f\colon X\to T$. 
Using \cite[Proposition 2.14]{Ishiki2021ultra}
or applying Theorem \ref{thm:Ishiki} to $A=\emptyset$, 
we can 
take  $d\in \ult(X; S)$. 
Then, Lemma \ref{lem:dormfult} implies that 
$d\lor \orddis_{T}[f]$ is a proper metric in $\ult(X; S)$. 
\end{proof}

Using Theorem 
\ref{thm:TUproper}, 
we obtain an extension theorem of proper metrics. 
\begin{thm}\label{thm:extproper}
Let $X$ be a
 $\sigma$-compact locally compact
Hausdorff 
  space, 
  and 
$A$  a 
non-empty
 closed  subset of $X$. 
If  
$d\in \met(A)$ is  a proper metric, 
then there exists a 
proper
metric $D\in \met(X)$ with 
$D|_{A^{2}}=d$. 
\end{thm}

\begin{proof}
Fix $p\in A$ and  define
a map  $f\colon A\to [0, \infty)$
 by
 $f(x)=d(p, x)$. Then $f$ is a
 continuous  proper function. 
 According to  Theorem \ref{thm:TUproper}, 
 there exists a continuous proper function 
 $F\colon X\to [0, \infty)$ with $F|_{A}=f$. 
 Due to   Hausdorff's extension theorem
 (Theorem \ref{thm:Hausdorff}), 
we can take a metric $e\in \met(X)$ such 
that $e|_{A^{2}}=d$. 
We define a map $D\colon X^{2}\to [0, \infty)$ by 
\[
D(x, y)=e(x, y)\lor \eucdis[F](x, y)
\]
Lemma \ref{lem:dormf} implies that 
the map $D$ is a proper metric  in $\met(X)$. 
We shall prove that  $D|_{A^{2}}=d$. 
If $x, y\in A$, we have 
$e(x, y)=d(x, y)$ and 
$\eucdis[F](x, y)=|F(x)-F(y)|=|d(x, p)-d(y, p)|$. 
The triangle inequality yields 
$|d(x, p)-d(y, p)|\le d(x, y)$. 
Thus, we obtain $\eucdis[F](x, y)\le d(x, y)$ for all $x, y\in A$. 
Therefore, by the definition of $D$, 
we have $D|_{A^{2}}=d$. 
 This completes the proof. 
\end{proof}

Let $X$ be a topological space. 
A subset $A$ of $X$ is said to be 
a \emph{retract} if there exists a 
continuous map $r\colon X\to A$ such that 
$r(a)=a$ for all $a\in A$. In this case, 
the continuous map $r$ is said to be 
a \emph{retraction}. 
A subset $A$ is said to be 
a \emph{proper retract} if there exists 
a retraction $r\colon X\to A$,  which is 
a proper map. 
For more discussion  of  proper retracts, 
we refer the readers to \cite{MR1909003}.

The  next lemma follows from the strong triangle inequality. 
\begin{lem}\label{lem:isosceles}
Let 
$X$ 
be a set, 
and
$d$ 
be an ultrametric on $X$. 
Then 
for all 
$x, y, z\in X$, 
the inequality 
$d(x, z)<d(y, z)$ 
implies 
$d(y, z)=d(x, y)$. 
\end{lem}

Let $(X, d)$ and $(Y, e)$ be 
metric spaces,
and $f\colon X\to Y$ be a map. 
We say that $f$ is \emph{metrically proper} if 
the inverse image $f^{-1}(A)$ is bounded in 
$(X, d)$
for every bounded subset $A$ of $Y$.

The proof of the following theorem is presented in 
\cite[Theorem 2.9]{brodskiy2007dimension}. 

\begin{thm}\label{thm:BDHM}
Let $(X, d)$ be an ultrametric space, 
and $A$ be a closed subset of $X$. 
Let $\tau\in (1, \infty)$. 
Then there exists a 
 $\tau^{2}$-Lipschitz
 retraction 
$r\colon X\to A$.  
 Moreover, 
if $A$ is unbounded, the retraction  $r$ associated with $A$ can be 
chosen to be  metrically proper. 
\end{thm}

By proving the existence of a proper ultrametric on 
an ultranormal  $\sigma$-compact 
locally compact Hausdorff  space
(Corollary \ref{cor:extpropermetult}), 
we show that a non-compact  closed subset 
of an ultranormal 
 $\sigma$-compact locally compact 
Hausdorff  
space  is not only just a retract, 
but also a proper retract.

\begin{thm}\label{thm:properretract}
Let $X$ be 
an ultranormal 
$\sigma$-compact locally compact
Hausdorff 
 space, 
 and 
$A$ a non-empty 
non-compact 
closed subset of $X$. 
Then 
$A$ is
a proper retract of $X$. 
\end{thm}
\begin{proof}
Using  Corollary \ref{cor:extpropermetult}, 
we can take 
a proper ultrametric  $d\in \ult(X; [0, \infty))$. 
Since $A$ is non-compact and $d$ is proper, 
it is unbounded in $(X, d)$. 
The latter part of Theorem \ref{thm:BDHM}
 implies that 
there exists a 
metrically proper retraction $r\colon X\to A$ 
with respect to  $d$. 
To prove that $r$ is proper, we take an arbitrary  compact subset $K$ of $A$. 
Since $K$ is bounded, and since $r$ is metrically proper, 
the inverse image $r^{-1}(K)$ is bounded and closed. 
Since $d$ is a proper metric, the set $r^{-1}(K)$ is compact, 
and hence $r$ is proper. 
This finishes the proof
of Theorem \ref{thm:properretract}. 
\end{proof}

Before proving the following corollary, 
notice that the composition of two proper maps is proper. 
\begin{cor}\label{cor:propermapext}
Let $X$ be an ultranormal 
$\sigma$-compact locally compact 
Hausdorff space, 
and $A$  a non-empty closed subset of $X$. 
If 
 $Y$ is  a non-compact metrizable space, 
then 
every
continuous  proper map $f\colon A\to Y$
can be extended into 
a
continuous 
proper map $F\colon X\to Y$. 
\end{cor}
\begin{proof}
We divide the proof into two cases. 

Case 1.  ($A$ is non-compact): 
Theorem \ref{thm:properretract}
guarantees the existence of a proper retraction  $r\colon X\to A$.  
Put  $F=f\circ r$. 
Then $F\colon X\to Y$ is a desired extension. 

Case 2.  ($A$ is compact): 
In this case, let  $Z$ be 
the countable discrete space. 
Fix $\omega\in A$. 
Put $Z=\{\, a_{i}\mid i\in \zz_{\ge 0}\, \}$, 
where $a_{*}\colon \zz_{\ge 0}\to Z$ is injective. 
Note that $X\times Z$ is an ultranormal 
non-compact $\sigma$-compact locally compact 
Hausdorff space. 
Put 
$C=A\times \{a_{0}\}\cup \{\omega\}\times Z$. 
Then $C$ is a non-compact 
closed  subset of $X\times Z$. 
Since $Y$ is non-compact, 
we can take a  countable closed discrete subset 
$\{\, b_{i}\mid i\in \zz_{\ge 1}\, \}$ of $Y$. 
We define a map $g\colon C\to Y$ by 
$g((x, a_{0}))=f(x)$ for all $x\in A$ and 
$g((\omega, a_{i}))=b_{i}$ for all $i\in \zz_{\ge 1}$. 
Then $g$ is continuous and  proper. 
Thus, using  Case 1, 
we can take  a continuous proper map 
$G\colon X\times Z\to Y$ such that $G|_{C}=g$. 
We define a map $F\colon X\to Y$ by 
$F(x)=G(x, a_{0})$. 
Then $F$ is a continuous proper map and 
satisfies $F|_{A}=f$. 
\end{proof}

\begin{prop}\label{prop:psimap}
Let $S$ be an unbounded characteristic subset of 
$[0, \infty)$, 
and $T$ a sporadic subset of $[0, \infty)$
with $T\yosub S$. 
Let $X$ be an ultranormal 
$\sigma$-compact locally compact
Hausdorff space. 
If  $d\in \ult(X; S)$, 
then there exists a 
metric $w\in \ult(X; T)$
such that $w(x, y)\le d(x, y)$
for all $x, y\in X$. 
Moreover, if $d$ is proper, so is $w$. 
\end{prop}
\begin{proof}
Take a real sequence  $\{a_{n}\}_{n\in \zz}$ such that 
$T=\{0\}\cup \{\, a_{n}\mid n\in \zz\, \}$,  
$\lim_{n\to \infty}a_{n}=\infty$, 
$\lim_{n\to -\infty}a_{n}=0$, and 
$a_{i}<a_{i+1}$ for all $i\in \zz$. 
We define a map $\psi\colon [0, \infty)\to [0, \infty)$ by 
\[
\psi(x)=
\begin{cases}
0 & \text{if $x=0$;}\\
a_{i} &\text{if $a_{i}\le x<a_{i+1}$.}
\end{cases}
\]
Put $w=\psi\circ d$. 
According to \cite[Lemma 2.2]{Ishiki2021ultra}, 
we observe that $w\in \ult(X; T)$. 
By the definition of $\psi$, 
we have $w(x, y)\le d(x, y)$ for all $x, y\in X$. 
This completes the first part of the proposition. 
To prove the latter part, 
assume that $d$ is proper and  take 
 $p\in X$ and  $r\in (0, \infty)$. 
 Put $\psi(r)=a_{i}$. 
 Then we have 
 $B(p, r; w)=B(p, a_{i}; w)\yosub B(p, a_{i+1}; d)$. 
 Since $d$ is proper, the set $B(p, r; w)$ is compact. 
 Thus $w$ is proper. 
\end{proof}

Theorem \ref{thm:properretract} provides 
 an ultrametric version  of Theorem \ref{thm:extproper}. 
\begin{thm}\label{thm:extproperult}
Let $S$ be an unbounded characteristic subset of 
  $[0, \infty)$. 
Let $X$ be 
an
ultranormal 
 $\sigma$-compact locally compact
Hausdorff 
  space, 
  and $A$  a 
  non-empty
 closed  subset of $X$. 
If 
$d\in \ult(A; S)$ is 
proper, 
then 
there exists a 
proper
metric $D\in \ult(X; S)$ such that 
$D|_{A^{2}}=d$. 
\end{thm}
\begin{proof}
The proof is similar to that of  Theorem \ref{thm:extproper}. 
Fix $p\in A$. 
Take a sporadic subset $T$ of $[0, \infty)$ 
with $T\yosub S$. 
Using Proposition \ref{prop:psimap}, 
we can take $w\in \ult(X; T)$ with 
$w(x, y)\le d(x, y)$ for all $x, y\in X$. 
We define  a map $f\colon A\to T$ by 
$f(x)=w(p, x)$. 
Then $f$ is a continuous proper function. 
According to  Corollary \ref{cor:propermapext}, 
we can take  a continuous proper function $F\colon X\to T$ 
such that $F|_{A}=f$. 
By 
Theorem \ref{thm:Ishiki}, 
 there exists a metric $e\in \ult(X; S)$
such that $e|_{A^2}=d$. 
We define a map $D\colon X^{2}\to S$ by 
\[
D(x, y)=e(x, y)\lor \orddis_{T}[F](x, y). 
\]
Lemma \ref{lem:dormfult} implies that
the map $D$ is 
a proper metric in $\ult(X; S)$. 
We shall prove $D|_{A^{2}}=d$. 
Take  $x, y\in A$. 
We may assume that $w(p, x)\le w(p, y)$. 
If $w(p, x)<w(p, y)$, 
Lemma \ref{lem:isosceles} 
yields
 $w(x, y)=w(p, y)$. 
Thus $\orddis_{T}[F](x, y)=w(x, y)\le d(x, y)$. 
If $w(p, x)=w(p, y)$, then, 
by the definition of $\orddis_{T}$, 
we have $\orddis_{T}[F](x, y)=0\le d(x, y)$. 
Thus, due to  $e|_{A^{2}}=d$ and the definition of $D$, 
we obtain  $D|_{A^2}=d$. 
This finishes the proof of 
Theorem \ref{thm:extproperult}. 
\end{proof}

\section{Proper metrics at large scales}

For  a topological  space 
$Y$, 
we denote by 
$\yopow{Y}$ 
the set of all non-empty   subsets  of $Y$. 
For topological spaces $X$ and $Y$, 
we say that a map 
$\phi\colon X\to \yopow{Y}$ 
is 
\emph{lower semi-continuous} 
if for every open subset 
$O$ of $Y$, 
the set 
$\{\, x\in X\mid \phi(x)\cap O\neq \emptyset\, \}$ 
is open in $X$. 
For a map $\phi\colon X\to \yopow{Y}$, 
a map $f\colon X\to Y$ is said to be 
a \emph{selection of $\phi$} if it is continuous and 
satisfies $f(x)\in \phi(x)$ for all $x\in X$.

The following proposition from  
E. Michael 
\cite[Proposition 1.4]{MR77107}
 states that the existence of a selection of 
a set-valued map is 
equivalent to the extension of 
a selection defined on a closed subset of the domain. 
\begin{prop}\label{prop:12equiv12}
Let $X$ and $Y$ be topological space. 
If $\mathcal{S}$ is  a subset of $\yopow{Y}$ 
containing all one-point subsets of $Y$, 
then the following statements are equivalent to each other: 
\begin{enumerate}
\item 
For all  lower semi-continuous map $\phi\colon X\to \mathcal{S}$, 
there exists a selection of $\phi$. 
\item If $A$ is a closed subset  of $X$, 
and $\phi\colon X\to \mathcal{S}$ is a lower semi-continuous map, 
and 
if   $f\colon X\to Y$ is  a selection  of the restricted  map 
$\phi|_{A}\colon A\to \mathcal{S}$, 
then there exists a map  $F\colon X\to Y$, 
 which is  
a selection of 
$\phi\colon X\to \mathcal{S}$ such that $F|_{A}=f$. 
\end{enumerate}
\end{prop}

Let 
$V$ 
be a Banach space. 
We  denote by 
$\conv(V)$ 
the set of all non-empty closed convex subsets  of $V$. 
The next theorem is known as 
Michael's selection theorem on paracompact spaces
(see \cite[Theorem 3.2$''$]{MR77107}): 
\begin{thm}\label{thm:mselection1}
Let $X$ be a paracompact space, 
and  $V$ a Banach space. 
If 
$\phi\colon X\to \conv(V)$ is 
a lower semi-continuous map, 
then there exists a selection 
of $\phi$. 
\end{thm}

For a topological space $Z$, 
we denote by 
$\clo(Z)$ 
the set of all non-empty closed  subsets  of $Z$. 
Recall that 
every ultranormal paracompact  space is $0$-dimensional 
(i.e., it has  covering dimension $0$).
The following theorem is 
known as 
the 
$0$-dimensional Michael selection theorem
(see \cite[Theorem 2]{MR1529282}): 
\begin{thm}\label{thm:mselection0}
Let $X$ be a $0$-dimensional paracompact space, 
$Z$ a completely metrizable space. 
If 
$\phi\colon X\to \clo(Z)$ is 
 a lower semi-continuous map, 
 then 
 there exists a selection of 
 $\phi$. 
\end{thm}
For every Banach space $V$ (resp.~completely metrizable space $Z$),
the set $\conv(V)$ (resp.~$\clo(Z)$) 
contains all one-point sets of $V$ (resp.~$Z$). 
Thus, we can apply
Proposition \ref{prop:12equiv12}
to Theorems
\ref{thm:mselection1} and \ref{thm:mselection0}, 
and then we obtain the next two theorems on 
extending selections:

\begin{thm}\label{thm:Michael}
Let $X$ be a paracompact space, 
and $A$  a closed subset of $X$. 
Let $V$ be a Banach space, 
and 
$\phi\colon X\to \conv(V)$ 
a lower semi-continuous map. 
If 
$f\colon A\to V$ 
is a selection of $\phi|_{A}\colon A\to \conv(V)$, 
then there exists a selection 
$F\colon X\to V$ of $\phi$
such that  $F|_A=f$. 
\end{thm}

\begin{thm}\label{thm:zeroMichael}
Let $X$ be an ultranormal paracompact space, 
and $A$  a closed subset of $X$. 
Let $Z$ be a completely metrizable space, 
and 
$\phi\colon X\to \clo(Z)$ 
 a lower semi-continuous map. 
If  
$f\colon A\to Z$ 
is a selection of $\phi|_{A}\colon A\to \clo(Z)$, 
then there exists a 
selection  
$F\colon X\to Z$  of $\phi$
such that  $F|_A=f$. 
\end{thm}

Propositions 
\ref{prop:lsc}
can be 
deduced from 
\cite[Theorem 0.48]{MR1659914} or 
\cite[Lemma 1.4.6]{MR977744}.  
The proof of Proposition \ref{prop:ultlsc}
is presented in 
\cite[Corollary 2.24]{Ishiki2021ultra}. 
The definition of ultra-normed modules can be found 
in  \cite{Ishiki2021ultra}.
\begin{prop}\label{prop:lsc}
Let $X$ be a topological space,  
and let $(V, \|*\|)$ be a Banach space. 
Let 
$H\colon X\to V$ 
be a continuous map and 
$r\in (0, \infty)$. 
Then the  map 
$\phi\colon  X\to \conv(V)$ 
defined by 
$\phi(x)=B(H(x), r; \|*\|)$ 
is lower semi-continuous. 
\end{prop}

\begin{prop}\label{prop:ultlsc}
Let $X$ be a topological space,  
$R$ be a commutative ring, 
and let $(V, h)$ be an ultra-normed $R$-module. 
Let 
$H\colon X\to V$ 
be a continuous map and 
$r\in (0, \infty)$. 
Then a map 
$\phi\colon X\to \clo(V)$ 
defined by 
$\phi(x)=B(H(x), r; h)$ 
is lower semi-continuous. 
\end{prop}

Let $(Z, h)$ be a metric space and  $\eta\in (0, \infty)$. 
A subset $E$ of $Z$ is said to be 
\emph{$\eta$-dense} in $(Z, h)$ if 
for all $x\in Z$, there exists $y\in E$ such that 
$h(x, y)\le \eta$.

\begin{thm}\label{thm:extdenseproper}
Let $\eta\in [0, \infty)$. 
Let $X$ be a $\sigma$-compact locally compact
Hausdorff  space, and $A$  a proper retract of $X$. 
If  $d\in \met(X)$ is  a proper metric, 
then 
there exists a proper metric 
$D\in \met(X)$ such that $D|_{A^{2}}=d$
and $A$ is $\eta$-dense in 
$(X, D)$. 
\end{thm}
\begin{proof}

We first  take a Banach space $(V, \|*\|)$
and an isometric embedding 
$l\colon (A, d|_{A^{2}})\to (V, \|*\|)$. 
For example, we can choose  $(V, \|*\|)$ as the space of 
all real-valued 
bounded continuous functions on $A$, and 
 $l\colon A\to V$ as the Kuratowski embedding defined by 
 $l(x)(y)=d(x, y)-d(\xi, y)$ for
 a  fixed point $\xi\in A$.

We  take a proper retraction $r\colon X\to A$
and 
define $\phi\colon X\to \conv(V)$ by 
$\phi(x)=B(l(r(x)), \eta; \|*\|)$. 
Applying Proposition \ref{prop:lsc} to $H=l\circ r$, 
we can assert that
the map  $\phi$ is lower semi-continuous. 
For all $a\in A$, 
the equality $r(a)=a$ implies that 
$l(a)\in \phi(a)$, 
namely, the map $l\colon A\to V$ is a 
selection of $\phi|_{A}\colon A\to \conv(V)$. 
Then 
Theorem \ref{thm:Michael} guarantees the 
existence of 
 a selection $L\colon X\to V$ of $\phi$ such that 
$L|_{A}=l$. 

Due to  Hausdorff's extension theorem
(Theorem \ref{thm:Hausdorff}), 
we can take $e\in \met(X)$ with 
$e|_{A^2}=d$. 
We define a map  $u\colon X^{2}\to [0, \infty)$ by  $u(x, y)=\min\{e(x, y), \eta\}$. 
Then $u\in\met(X)$. 
We also define a map $v\colon X^{2}\to [0, \infty)$ by  
$v(x, y)=\|L(x)-L(y)\|\lor u(x, y)$. 
Since  $u\in \met(X)$ and 
$L$ is continuous,  we have $v\in \met(X)$. 
Using  $L|_{A}=l$, 
we obtain  $\|L(a)-L(b)\|=d(a, b)$ for all $a, b\in A$. 
Then, from $u(x, y)\le e(x, y)$ for all $x, y\in X$, and  $e|_{A^{2}}=d$, 
it follows that
$v|_{A^{2}}=d$. 

Next we 
fix $p\in A$ (note that $A\neq \emptyset$). 
We define a continuous proper function 
$f\colon A\to [0, \infty)$ by 
$f(x)=d(p, x)$, 
and 
define $F=f\circ r$. 
Then $F\colon X\to [0, \infty)$ 
is a continuous  proper
function with $F|_{A}=f$.  
We also define a metric $D$ on $X$ by 
$D(x, y)=v(x, y)\lor \eucdis[F](x, y)$. 

Lemma \ref{lem:dormf} implies that 
$D$ is a proper metric in  $\met(X)$. 
In  a similar way to the proof of Theorem \ref{thm:extproper}, 
we obtain $D|_{A^{2}}=d$. 

We now show that 
$A$ is $\eta$-dense in $(X, D)$. 
Take an arbitrary point  $x\in X$. 
The relations $L(r(x))=l(r(x))$ and $L(x)\in \phi(x)$ 
yield
\begin{align}\label{al:llreta}
\|L(x)-L(r(x))\|\le \eta.
\end{align} 
From   \eqref{al:llreta}, the inequality $u(x, r(x))\le \eta$, 
and the definition of $v$, 
it follows that
\begin{align}\label{al:vrx}
v(x, r(x))\le \eta.
\end{align}
Since $r$ is a retraction, 
we have $r(r(x))=r(x)$. 
Thus 
$\eucdis[F](x, r(x))=|F(x)-F(r(x))|=|f(r(x))-f(r(x))|=0$, 
and hence 
\begin{align}\label{al:eff}
\eucdis[F](x, r(x))=0. 
\end{align}
Therefore, 
by \eqref{al:vrx}, \eqref{al:eff},
and the definition of $D$,  we conclude that \begin{align}\label{al:dxx}
D(x, r(x))\le \eta.
\end{align} 
Since $r(x)\in A$, 
and $x\in X$ is arbitrary, 
the inequality \eqref{al:dxx}
proves that $A$ is $\eta$-dense in 
$(X, D)$. 
This completes the proof of 
Theorem \ref{thm:extdenseproper}. 
\end{proof}

The proof of Theorem \ref{thm:extultproper} is 
analogous with Theorems \ref{thm:extdenseproper}.

\begin{thm}\label{thm:extultproper}
Let $\eta\in (0, \infty)$,
and 
$S$  an unbounded characteristic subset of 
$[0, \infty)$. 
Let $X$ be an
ultranormal 
$\sigma$-compact locally compact
Hausdorff 
 space, 
 and 
 $A$  a non-empty 
non-compact 
closed subset of $X$. 
If  $d\in \ult(A; S)$ is 
proper, 
then 
there
exists a 
proper metric 
$D\in \ult(X; S)$
such that 
$D|_{A^{2}}=d$ and 
$A$ is $\eta$-dense in 
$(X, D)$. 
\end{thm}

\begin{proof} 
We put  $R=\zz/2\zz$. 
However, 
as long as $R$ is an integral domain, 
the choice of $R$ does not affect the 
proof of the theorem. 
We first verify that there exists  an isometric embedding 
$(A, d|_{A^{2}})$ into 
a complete ultra-normed $R$-module. 
Let $(Y, m)$ be the completion of $(A, d|_{A^{2}})$. 
Since the set $d(A^{2})$ is invariant under the completion (see (12) in \cite[Theorem 1.6]{MR3782290}), 
we have $m(Y^{2})=d(A^{2})$, and hence
 $m\in \ult(Y; S)$. 
According to \cite[Theorem 1.1]{Ishiki2021ultra}, 
we can take 
a complete ultra-normed $R$-module
 $(V,  h)$
 with $h\in \ult(V; S)$ and an isometric embedding 
$J\colon (Y, m)\to (V, h)$. 
We put 
$l=J|_{A}\colon A\to V$.

Theorem \ref{thm:properretract} enables us to 
 take a proper retraction $r\colon X\to A$. 
Since $S$ is characteristic, we can also  take 
$\theta\in S\setminus \{0\}$ with $\theta\le \eta$. 
We define a map 
$\phi\colon X\to \clo(V)$ by 
$\phi(x)=B(l(r(x)), \theta; h)$. 
Applying   Proposition  \ref{prop:ultlsc} to $H=l\circ r$, 
we notice that the map $\phi$ is lower semi-continuous. 
For all $a\in A$, the equality 
$r(a)=a$ implies that
$l(a)\in \phi(a)$, 
namely, the map $l\colon A\to V$ is a 
selection of $\phi|_{A}\colon A\to \clo(V)$. 
Using Theorem \ref{thm:zeroMichael}, 
there exists a selection  $L\colon X\to V$ of $\phi$
such that 
$L|_{A}=l$. 

Due to  Theorem \ref{thm:Ishiki}, 
we can take $e\in \ult(X; S)$ such that 
$e|_{A^2}=d$. 
Put $u(x, y)=\min\{e(x, y), \theta\}$. 
Then $u\in \ult(X; S)$. 
We define a map $v\colon X^{2}\to [0, \infty)$ by  
$v(x, y)=h(L(x), L(y))\lor u(x, y)$. 
From $u\in \ult(X; S)$
and the continuity of $L$, 
it follows that 
$v\in \ult(X; S)$. 
Using   $L|_{A}=l$, 
we have 
$h(L(a), L(b))=d(a, b)$ for all $a, b\in A$. 
Then, by  $u(x, y)\le e(x, y)$ for all $x, y\in X$, 
and by $e|_{A^{2}}=d$, 
we obtain  $v|_{A^{2}}=d$. 

Next we 
fix $p\in A$ and 
take a sporadic subset $T$ of $[0, \infty)$ with 
$T\yosub S$. 
Due to Proposition \ref{prop:psimap}, 
there exists $w\in \ult(A; T)$ with 
$w(a, b)\le d(a, b)$ for all $a, b\in A$. 
We define 
a continuous proper  function  $f\colon A\to T$ by 
$f(x)=w(p, x)$, 
and define
a map  $F=f\circ r$. 
Then $F\colon X\to T$ 
is a continuous  proper
function
with $F|_{A}=f$.  
We also define a metric $D$ on $X$ by 
$D(x, y)=v(x, y)\lor \orddis_{T} [F](x, y)$. 

Lemma \ref{lem:dormfult} implies that 
$D$ is a proper ultrametric in  $\ult(X; S)$. 
In  a similar way  to 
the proof of Theorem \ref{thm:extproperult}, 
we obtain $D|_{A^{2}}=d$. 

We now show that 
$A$ is $\eta$-dense in $(X, D)$. 
Take an  arbitrary point  $x\in X$. 
The relations $L(r(x))=l(r(x))$ and   $L(x)\in \phi(x)$
yield
\begin{align}\label{al:hll}
h(L(x), L(r(x)))\le \theta. 
\end{align}\label{al:utheta}
From \eqref{al:hll}, 
the inequalities $u(x, r(x))\le \theta$ and
$\theta\le \eta$, and the definition of $v$, 
it follows that 
\begin{align}\label{al:ultvxrx}
v(x, r(x))\le  \eta.  
\end{align}
Since $r$ is a retraction, 
we have $r(r(x))=r(x)$. 
Then 
$\orddis_{T}[F](x, r(x))=\orddis_{T}(F(x), F(r(x)))
=\orddis_{T}(f(r(x)), f(r(x)))=0$, 
and hence 
\begin{align}\label{al:mtxrx}
\orddis_{T}[F](x, r(x))=0. 
\end{align}
Therefore, by 
\eqref{al:ultvxrx}, \eqref{al:mtxrx}, 
and the definition of $D$,  we conclude that 
\begin{align}\label{al:ddxrx}
D(x, r(x))\le \eta. 
\end{align}
Since $r(x)\in A$
and $x\in X$ is arbitrary, 
the inequality \eqref{al:ddxrx} 
proves that $A$ is $\eta$-dense in 
$(X, D)$. This 
completes the proof
of 
Theorem \ref{thm:extultproper}. 
\end{proof}

\begin{ac}
The author would like to thank 
the referee for helpful comments and suggestions. 
\end{ac}


\bibliographystyle{amsplain}
\bibliography{bibtex/proper.bib}

\providecommand{\bysame}{\leavevmode\hbox to3em{\hrulefill}\thinspace}
\providecommand{\MR}{\relax\ifhmode\unskip\space\fi MR }
\providecommand{\MRhref}[2]{%
  \href{http://www.ams.org/mathscinet-getitem?mr=#1}{#2}
}
\providecommand{\href}[2]{#2}
\begin{thebibliography}{10}

\bibitem{Arens1952}
R.~Arens, \emph{Extension of functions on fully normal spaces}, Pacific. J.
  Math. \textbf{2} (1952), 11--22.

\bibitem{Bing1947}
R.~H. Bing, \emph{Extending a metric}, Duke. Math. J. \textbf{14} (1947),
  no.~3, 511--519.

\bibitem{brodskiy2007dimension}
N.~Brodskiy, J.~Dydak, J.~Higes, and A.~Mitra, \emph{Dimension zero at all
  scales}, Topology Appl. \textbf{154} (2007), no.~14, 2729--2740. \MR{2340955}

\bibitem{MR3782290}
A.~B. Comicheo and K.~Shamseddine, \emph{Summary on non-{A}rchimedean valued
  fields}, Advances in {U}ltrametric {A}nalysis (A.~Escassut, C.~Perez-Garcia,
  and K.~Shamseddibe, eds.), Contemp. Math., vol. 704, Amer. Math. Soc.,
  Providence, RI, 2018, pp.~1--36. \MR{3782290}

\bibitem{MR80905}
J.~de~Groot, \emph{Non-{A}rchimedean metrics in topology}, Proc. Amer. Math.
  Soc. \textbf{7} (1956), 948--953. \MR{80905}

\bibitem{MR2435142}
C.~Delhomm\'{e}, C.~Laflamme, M.~Pouzet, and N.~Sauer, \emph{Indivisible
  ultrametric spaces}, Topology Appl. \textbf{155} (2008), no.~14, 1462--1478.
  \MR{2435142}

\bibitem{MR2854677}
D.~Dordovskyi, O.~Dovgoshey, and E.~Petrov, \emph{Diameter and diametrical
  pairs of points in ultrametric spaces}, $p$-Adic Numbers Ultrametric Anal.
  Appl. \textbf{3} (2011), no.~4, 253--262. \MR{2854677}

\bibitem{MR3090172}
O.~Dovgoshey, O.~Martio, and M.~Vuorinen, \emph{Metrization of weighted
  graphs}, Ann. Comb. \textbf{17} (2013), no.~3, 455--476. \MR{3090172}

\bibitem{MR3135687}
O.~Dovgoshey and E.~A. Petrov, \emph{A subdominant pseudoultrametric on
  graphs}, Mat. Sb. \textbf{204} (2013), no.~8, 51--72. \MR{3135687}

\bibitem{Du1951}
J.~Dugundji, \emph{An extension of {T}ietze's theorem}, Pacific J. Math.
  \textbf{1} (1951), 353--367.

\bibitem{ellis1970extending}
R.~L. Ellis, \emph{Extending continuous functions on zero-dimensional spaces},
  Math. Ann. \textbf{186} (1970), no.~2, 114--122.

\bibitem{MR1346246}
M.~Frantz, \emph{Controlling {T}ietze-{U}rysohn extensions}, Pacific J. Math.
  \textbf{169} (1995), no.~1, 53--73. \MR{1346246}

\bibitem{Ha1930}
F.~Hausdorff, \emph{{E}rweiterung einer {H}om{\"o}omorphie}, Fund. Math.
  \textbf{16} (1930), 353--360.

\bibitem{Hu2010}
M.~Hu\v{s}ek, \emph{Extension of mappings and pseudometrics}, Extracta Math.
  \textbf{25} (2010), no.~3, 277--308.

\bibitem{Ishiki2021ultra}
Y.~Ishiki, \emph{An embedding, an extension, and an interpolation of
  ultrametrics}, $p$-Adic Numbers Ultrametric Anal. Appl. \textbf{13} (2021),
  no.~2, 117--147.

\bibitem{MR77107}
E.~Michael, \emph{Continuous selections. {I}}, Ann. of Math. (2) \textbf{63}
  (1956), 361--382. \MR{77107}

\bibitem{MR1529282}
\bysame, \emph{Selected selection theorems}, Amer. Math. Monthly \textbf{63}
  (1956), no.~4, 233--238. \MR{1529282}

\bibitem{MR1909003}
\bysame, \emph{Closed retracts and perfect retracts}, Topology Appl.
  \textbf{121} (2002), no.~3, 451--468. \MR{1909003}

\bibitem{MR3728284}
J.~R. Munkres, \emph{Topology}, 2nd ed., Pearson modern classic, Pearson, New
  York, 2018, Originally published in 2000 [2018 reissue].

\bibitem{MR1659914}
D.~Repov\v{s} and P.~V. Semenov, \emph{{C}ontinuous {S}elections of
  {M}ultivalued {M}appings}, Mathematics and its Applications, vol. 455, Kluwer
  Academic Publishers, Dordrecht, 1998. \MR{1659914}

\bibitem{MR321026}
H.~Toru\'{n}czyk, \emph{A short proof of {H}ausdorff's theorem on extending
  metrics}, Fund. Math. \textbf{77} (1972), no.~2, 191--193. \MR{321026}

\bibitem{MR977744}
J.~van Mill, \emph{{I}nfinite-{D}imensional {T}opology: {P}rerequisites and
  {I}ntroduction}, North-Holland Mathematical Library, vol.~43, North-Holland
  Publishing Co., Amsterdam, 1988. \MR{977744}

\bibitem{MR1967003}
K.~Yamazaki, \emph{Controlling extensions of functions and {$C$}-embedding},
  Topology Proc. \textbf{26} (2001/02), no.~1, 323--341. \MR{1967003}

\end{thebibliography}

\end{document}